\documentclass[11pt]{amsart} 
\usepackage{amssymb,amsmath,graphicx}

\newtheorem{theorem}{Theorem}[section]
\newtheorem{corollary}[theorem]{Corollary}

\title[Ehrhart-Macdonald reciprocity extended]{Ehrhart-Macdonald reciprocity extended} 

\author{Matthias Beck}
\address{Department of Mathematics\\
         San Francisco State University\\
         San Francisco, CA 94132\\
         USA}
\email{beck@math.sfsu.edu}
\urladdr{http://math.sfsu.edu/beck}

\author{Richard Ehrenborg}
\address{Department of Mathematics\\
         University of Kentucky\\
         Lexington, KY 40506-0027 \\
         USA}
\email{jrge@ms.uky.edu}
\urladdr{http://www.ms.uky.edu/\~{}jrge/}

\thanks{We thank Margaret Readdy for helpful comments.
Parts of this paper was written while the first author was at SUNY Binghamton and the Mathematical Sciences Research Institute, Berkeley; he thanks both institutions for their hospitality. The second author was partially supported by National Science Foundation grant 0200624 and by a University of Kentucky College of Arts and Sciences 2004 summer grant.}


\date{April 7, 2005}
\begin{document}
\renewcommand{\P}{{\mathcal P}}

\begin{abstract}
For a convex polytope $\P$ with rational vertices, we count the number of integer points in integral dilates of $\P$ and its interior. The Ehrhart-Macdonald reciprocity law gives an intimate relation between these two counting functions. A similar counting function and reciprocity law exists for the sum of all solid angles at integer points in dilates of $\P$. We derive a unifying generalization of these reciprocity theorems which follows in a natural way from Brion's Theorem on conic decompositions of polytopes. 
\end{abstract}

\maketitle
%\tableofcontents

%%%%%%%%%%%%%%%%%%%%%%%%%%%%%%%%%%%%%%%%%%%%%%%%%%%%%%%%%%%%%%%%%%%%%%%%%%%
\newcommand{\F}{{\mathcal F}}
\newcommand{\K}{{\mathcal K}}
\newcommand{\R}{{\mathbb R}}
\newcommand{\Z}{{\mathbb Z}}
\newcommand{\m}{{\bf m}}
\newcommand{\x}{{\bf x}}

\section{Introduction}
\setcounter{equation}{0}

The origin of this work lies in two beautiful \emph{reciprocity
theorems} of Ehrhart and Macdonald. One attaches to a convex rational
polytope $\P$ two counting functions:\ the number of integral points
in $t \P$ and the sum of the solid angles in $t \P$, as functions of
a positive integer $t$. These functions are quasi-polynomials in $t$
which have a life beyond the positive integers. Namely, when evaluated
at negative integers, they give the respective counting function for
the \emph{interior} of $\P$. Our goal is to unify and generalize the
integral-point and solid-angle counting functions, using the
philosophy of valuations. Our main tool is Brion's Theorem on conic
decompositions of polytopes.

A \emph{convex polyhedron} $\P$ is the intersection of (open or closed) half-spaces in $\R^d$. We say $\P$ is a \emph{convex polytope} if it is bounded, and $\P$ is a \emph{cone} if each hyperplane bounding the half-spaces that define $\P$ contains the origin.
To each integer point $\m \in \Z^d$ in a convex polyhedron $\P$ we attach the \emph{tangent cone} $\K_\P (\m)$, defined as follows. Denote the \emph{characteristic function} of a set $X$ by $1_X(x)$, that is, it is 1 for $x \in X$ and 0 otherwise.
Let $f_\P (\m,\x) = \lim_{ \epsilon \to 0^+ } 1_\P ( \m + \epsilon \x )$. Then the tangent cone at $\m$ is
$$
\K_\P (\m) = \left\{ \x \in \R^d : \ f_\P ( \m, \x ) = 1 \right\} .
$$
In other words, $\K_\P (\m)$ equals
\begin{itemize}
\item $\R^d$ if $\m$ is in the interior of $\P$;
\item $\emptyset$ if $\m$ is not in the closure of $\P$;
\item the intersection of the halfspaces defining a face $\F$ if $\m$ is on the face $\F$ of $\P$, shifted to the origin.
\end{itemize}

\renewcommand{\L}{{\mathcal L}}

A \emph{lattice of sets} is a collection $\L$ of sets, partially ordered by set inclusion, such that for any two sets $A, B \in \L$ 
  \[ A \cap B \in \L \qquad \text{ and } \qquad A \cup B \in \L \ . \]
A function $v$ on $\L$ with values in some abelian group is a \emph{valuation} (often called a \emph{finitely additive measure}) if for any $A, B \in \L$ 
  \begin{equation}\label{inclexcl} v ( A \cup B ) = v(A) + v(B) - v ( A \cap B ) \ . \end{equation} 
We consider the lattice $\L$ generated by (open, closed, half-open) cones. 
We define $v_\P (\m)$ as the valuation $v$ evaluated at the tangent cone $\K_\P (\m)$. This is a valuation when viewed as a function of $\P$. However, we will mostly be interested in $v_\P (\m)$ as a function in $\m$.

The counting function we study is
\[
V_\P (t) := \sum_{ \m \in \Z^d } v_{t \P} (\m) \ ,
\]
initially defined for positive integers $t$.
This function only makes sense if $\P$ is bounded, that is, $\P$ is a convex polytope. In order to state our main theorem, recall that a \emph{quasi-polynomial} $Q$ of \emph{degree} $d$ as an expression of the form
  \[ Q(t) = c_{d}(t) \, t^{d} + \dots + c_{1}(t) \, t + c_{0}(t) \ , \]
where $c_0,c_1,\dots,c_d$ are periodic functions of $t$ and $c_d \not\equiv 0$. 
The least common multiple of the periods of the $c_j$ is called the \emph{period} of $Q$.
The \emph{denominator} of the polytope $\P$ is the least common multiple of the denominators of the vertices of $\P$.
We denote the relative interior of $\P$ by $\P^\circ$. 

\begin{theorem}\label{T:mainthm}
For any closed convex rational $d$-polytope, $V_\P (t)$ is a quasi-polynomial in $t$ of degree~$d$. The period of $V_\P (t)$ divides the denominator of $\P$. In particular, if $\P$ has integer vertices then $ V_\P $ is a polynomial. Its evaluation at negative integers gives 
\begin{equation}\label{E:mainrec}
V_\P (-t) = (-1)^d \, V_{ - \P^\circ } (t) \ .
\end{equation}
\end{theorem}

Our motivation for studying the function $V_\P (t)$ comes from special cases of valuations. The first example is
\[
v(A) := \begin{cases}
1 & \text{ if } 0 \in A , \\
0 & \text{ if } 0 \not\in A .
\end{cases}
\]
In this case $V_\P (t)$ counts the number of integer points in $ t \P $ and Theorem~\ref{T:mainthm} specializes to the following fundamental theorem due to Ehrhart~\cite{ehrhartpolynomial,ehrhart1}. Note that $V_\P (t) = V_{ -\P } (t)$.

\begin{corollary}[Ehrhart]\label{C:ehrhart} 
If $P$ is a closed convex rational $d$-polytope, then $ L_\P (t) := \# \left( t \P \cap \Z^d \right)  $ is a quasi-polynomial of degree $d$, having period that divides the denominator of $\P$. In particular, if $\P$ has integer vertices, then $ L_\P $ is a polynomial. The evaluation of $L_\P$ at negative integers gives
\begin{equation}\label{E:ehrrec}
  L_\P ( -t ) = (-1)^d \, L_{\P^\circ} ( t )  \ . 
\end{equation}
\end{corollary} 
The \emph{reciprocity law} (\ref{E:ehrrec}), conjectured and partially proved by Ehrhart, was in its full generality first proved by Macdonald~\cite{macdonald}.

The second example of a valuation is 
\[
v(A) := \text{ solid angle of $A$ at $0$ } ,
\]
that is, the ratio of the volumes of $A \cap B$ and $B$, where $B$ is a sufficiently small $d$-ball centered at the origin.
Now $V_\P (t)$ counts the integer points in $ t \P $ weighted by their solid angle in $t \P$. Theorem~\ref{T:mainthm} specializes to the following fundamental theorem due to Macdonald~\cite{macdonald}. Note again that $V_\P (t) = V_{ -\P } (t)$.

\begin{corollary}[Macdonald]\label{C:macdonald} 
Suppose $\P$ is a convex rational $d$-po\-lytope. Let $A_\P (t)$ be the sum of the solid angles of $t \P$ at all integer points. Then $A_\P$ is a quasi-polynomial of degree $d$, whose period divides the least common multiple of the denominators of the vertices of $\P$. The evaluation of $A_\P$ at negative integers gives
\begin{equation}\label{E:macrec}
  A_\P (-t) = (-1)^d \, A_\P (t) \ .
\end{equation}
In particular, if $\P$ has integer vertices, then $ A_\P $ is a polynomial which is either even or odd. 
\end{corollary} 
It is not hard to conclude the first half of Macdonald's Theorem from Ehrhart's Theorem (Corollary~\ref{C:ehrhart}), since all integer points in a given face of $\P$ have the same solid angle. The nontrivial part is equation (\ref{E:macrec}). 

\renewcommand{\v}{{\bf v}}
\newcommand{\D}{{\mathcal D}}

A third example of a valuation is as follows. Fix a vector $\v \in \R^d$, and let
\[
v(A) :=
\begin{cases}
1 & \text{ if } \v \in A , \\
0 & \text{ if } \v \not\in A .
\end{cases}
\]
This time $V_\P (t)$ counts those integer points $\m$ that after a small step in direction $\v$ will be in $t \P$.
Note that, in particular, $V_\P (t) \not= V_{ -\P } (t)$. Theorem~\ref{T:mainthm} specializes to the following statement. Let $\D (\P,\v) = \left\{ \m \in \Z^d : \ \m + \epsilon \v \in \P \text{ for small enough } \epsilon > 0 \right\}$ and $D_{ \P,\v } (t) = \# \D (t\P,\v)$.

\begin{corollary}
Suppose $\P$ is a convex rational $d$-po\-lytope and $\v$ is a fixed vector in $\R^d$. Then $D_{ \P,\v } (t)$ is a quasi-polynomial of degree $d$. The period of $D_{ \P,\v } (t)$ divides the least common multiple of the denominators of the vertices of $\P$. In particular, if $\P$ has integer vertices, then $D_{ \P,\v }$ is a polynomial. The evaluation of $D_{ \P,\v }$ at negative integers gives
\[
  D_{ \P,\v } (-t) = (-1)^d \, D_{ -\P^\circ, \v } (t)  \ .
\]
\end{corollary} 

We conceived Theorem~\ref{T:mainthm} when we realized that Ehrhart's and Macdonald's theorems (Corollaries~\ref{C:ehrhart} and~\ref{C:macdonald}) follow very naturally from a theorem of Brion~\cite{brion}, which in itself can be thought of as a consequence of a theorem due to Brianchon~\cite{brianchon} and Gram~\cite{gram}. To state these results, we need to introduce a few notions of polyhedral geometry in the next section. 
By no means do we claim that our approach is the easiest to Corollaries~\ref{C:ehrhart} and~\ref{C:macdonald}; our goal is merely to show the intimate connections of these results. Another such connection is given by Barvinok's result~\cite{barvinokalgorithm}, which uses Brion's Theorem to conclude that in fixed dimension, the generating function of an Ehrhart quasi-polynomial is polynomial-time computable.

%%%%%%%%%%%%%%%%%%%%%%%%%%%%%%%%%%%%%%%%%%%%%%%%%

\section{Conic decompositions of polytopes}\label{S:conicdecomp}
\setcounter{equation}{0}
\renewcommand{\a}{{\bf a}}

Given a $d$-dimensional closed convex polyhedron $ \P \subset \R^{d} $, the hyperplane $H = \{ \x \in \R^d : \ \a \cdot \x = b \}$ is a \emph{supporting hyperplane} of $ \P $ if $\P \cap H \not= \emptyset$ and $\P$ lies entirely on one side of $H$, that is, $ \P \subseteq \{ \x \in \R^d : \ \a \cdot \x \leq b \}$ or $ \P \subseteq \{ \x \in \R^d : \ \a \cdot \x \geq b \}$. 
A \emph{face} of $ \P $ is a set of the form $ \P \cap H $, where $H$ is a supporting hyperplane of $\P$. Note that both $ \P $ itself and the empty set are faces of $ \P $. The $(d-1)$-dimensional faces are called \emph{facets}, the 1-dimensional faces \emph{edges}, and the 0-dimensional faces \emph{vertices} of $\P$. 

For a vertex $\v$ of $\P$, we define the \emph{affine vertex cone} $\P_\v$ as the shifted cone $\K_\P(\v) + \v$.

\newcommand{\C}{{\mathbb C}}
\newcommand{\z}{{\bf z}}

The \emph{integer generating function} $\sigma_\P$ of a convex polyhedron $\P$ is 
  \[ \sigma_\P(\z) = \sigma_\P(z_1,\dots,z_d) = \sum_{ (m_1,\dots,m_d) \in \P \cap \Z^d } z_1^{m_1} \cdots z_d^{m_d} = \sum_{ \m \in \P \cap \Z^d } \z^\m  \ . \] 
Here we view $\z$ as a $d$-dimensional complex variable. 
%However, one can also view $z \in \R$.
%The sum on the right-hand side does not always converge; however, if $\P$ is a cone, the sum converges absolutely for all $\z$ in some open subset of $\C^d$.) 

A polyhedron is \emph{rational} if each of its defining halfspaces can be described as $\{ \x \in \R^d : \ \a \cdot \x \leq b \}$ for some $\a \in \Z^d$ and $b \in \Z$. If $\P$ is a rational polyhedron then $\sigma_\P(\z)$ can be written as a rational function in $z_1,\dots,z_d$ (see, for example, \cite{stanleyec1}). 
%Below we will derive the precise form of this rational function for the case that $\K$ is simple. 
The following theorem can be derived, for example, from the Brianchon-Gram Theorem~\cite{brianchon,gram}; this is shown in~\cite{conicdecomp}. The original proof of Brion is in~\cite{brion}; other, more elementary proofs can be found in~\cite{ishida,lawrence}. 

\begin{theorem}[Brion]\label{T:brion}
For a convex rational polytope $\P \subset \R^d$, represent the generating functions of the vertex cones $\P_\v$ as rational functions. Then 
  \[ \sum_{ \v \text{ \rm vertex of } \P } \sigma_{\P_\v} (\z) \ = \ \sigma_\P (\z) \ . \] 
\end{theorem} 

We need to extend Brion's Theorem to our more general setting, albeit for our purposes we only need a version for simplices. Our extension concerns the generating function for the valuation $v_\P$, that is,
\[
   \beta_\P (\z) = \sum_{ \m \in \Z^d } v_\P (\m) \, \z^\m \ .
\]
This is a rational function because $v_\P(\m)$ is constant on faces.

\renewcommand{\S}{{\mathcal S}}

\begin{theorem}\label{T:brionext}
Suppose $\S$ is a rational simplex and $v$ is a valuation. Then, as rational functions,
  \[ \sum_{ \v \text{ \rm vertex of } \S } \beta_{ \S_\v } (\z) \ = \ \beta_\S (\z) \ . \] 
\end{theorem} 

\newcommand{\W}{{\mathcal W}}
\newcommand\ch{\operatorname{ch}}
\newcommand{\G}{{\mathcal G}}

In what follows, we adjust a valuation-based proof of Brion's Theorem given in~\cite{conicdecomp} to our more general setting.

\begin{proof}[Proof of Theorem~\ref{T:brionext}]
Suppose $\S$ is a rational $d$-simplex. The hyperplanes bounding $\S$ divide $\R^d$ into regions. These regions are in one-to-one correspondence with the non-empty faces of $\S$, namely the closure of each region touches a unique maximal face of $\S$. Denote the region corresponding to the face $\F$ by $\ch(\F)$. We orient the hyperplanes in such a way that they point towards $\S$. Thus $\ch(\S) = \S$ is closed, the region corresponding to a vertex is open, and all the other regions are half open. 

Let $\W$ be the union of regions corresponding to subfaces of a face $\F \subseteq \S$, that is,
  \[ \W = \bigcup_{ \G \subseteq \F } \ch(\G) \ . \] 
The set $\W$ is a convex polyhedron that contains a line, except when $\F$ is a vertex. 
Because $v$ is a valuation, the generating function $\beta_\W (\z)$ satisfies,
  \[ \beta_\W (\z) = \sum_{ \G \subseteq \F } \beta_{\ch(\G)} (\z) \] 
Now unless $\F$ is a vertex, $\beta_\W (\z)$ is zero: There is an integer vector $\m$ such that $\W + \m = \W$; but 
$1 - \z^\m$ is not a zero divisor in the ring of rational functions in $\z$. 
In summary, we have 
\[
  \sum_{ \G \subseteq \F } \beta_{\ch(\G)} (\z) = 
  \begin{cases}
    \beta_{\ch(\F)} (\z) & \text{ if $\F$ is a vertex, } \\
    0 & \text{ otherwise. } 
  \end{cases}
\]
By inclusion-exclusion,
  \[ \beta_{\ch(\F)} (\z) = (-1)^{\dim \F} \! \! \sum_{\v\text{ \rm vertex of }\F} \! \beta_{\ch(\v)} (\z) \ . \] 
Again using that $v$ is a valuation, we have
  \begin{align*} \sum_{\v\text{ \rm vertex of }\S} \! \beta_{\P_\v} (\z) 
    &= (d+1) \beta_\S (\z) \, + \! \sum_{\F\text{ \rm facet of }\S} \! \beta_{\ch(\F)} (\z) \\ 
    &= \beta_\S (\z) + d (-1)^d \! \sum_{\v\text{ \rm vertex of }\S} \! \beta_{\ch(\v)} (\z)  + \! \sum_{\F\text{ \rm facet of }\S} \! (-1)^{d-1} \! \sum_{\v\text{ \rm vertex of }\F} \! \beta_{\ch(\F)} (\z) \\ 
    &= \beta_\S (\z) + \! \sum_{\v\text{ \rm vertex of }\S} \! \left( d (-1)^d + d (-1)^{d-1} \right) \beta_{\ch(\F)} (\z) \\ 
    &= \beta_\S (\z) \ . \end{align*} 
This completes the proof.
\end{proof}

Theorem~\ref{T:brionext} can be extended to general convex polytopes using the valuation theorems of the next section. Since we do not need the general case to prove our main theorem, we leave the proof to the reader.

\begin{theorem}
Suppose $\P$ is a convex rational polytope and $v$ is a valuation. As rational functions, we have the identity
  \[ \sum_{ \v \text{ \rm vertex of } \P } \beta_{ \P_\v } (\z) \ = \ \beta_\P (\z) \ . \] 
\end{theorem} 

%We will prove Theorem~\ref{T:brionext} for \emph{simplices}, i.e., $d$-polytopes with (a minimal number of) $d+1$ vertices. Once this is established, the theorems should follow by a triangulation argument. This last step has to be taken with a grain of salt---it is not entirely trivial (contrary to some statements in the literature). To this extend, we first introduce a valuation theorem of Volland, in the following section. Finally, we will deduce Theorem~\ref{T:mainthm} from Theorem~\ref{T:brionext}.

%%%%%%%%%%%%%%%%%%%%%%%%%%%%%%%%%%%%%%%%%%%%%%%%%%%%%%%%%%%%%%%%%%%%%%%%%%%%%%
\newcommand{\T}{{\mathcal T}}
\newcommand{\Rat}{\mathrm{Rat}}

\section{Valuations on the triangulation lattice} 
\setcounter{equation}{0}

A subset $\G$ of a lattice $\L$ is a \emph{generating set} if $\G$ is closed under intersection and every element of $\L$ is the union of elements of $\G$. The following \emph{extension theorem} is due to Groemer~\cite{groemer} and generalizes a valuation theorem for polyhedra due to Volland~\cite{volland}.

\begin{theorem}[Groemer]\label{T:groemer} 
Suppose $\G$ generates the lattice $\L$, and $v$ is a function on $\G$ that satisfies 
\[
  v \left( \bigcup_{ i=1 }^{n} A_{i} \right) = \sum_{ 1 \le i \le n } v \left( A_{i} \right) - \sum_{ 1 \le i<j \le n } v \left( A_{i} \cap A_{j} \right) + \sum_{ 1 \le i<j<k \le n } v \left( A_{i} \cap A_{j} \cap A_{k} \right) - \cdots
\]
whenever $A_1, A_2, \dots, A_n, A_1 \cup A_2 \cup \dots \cup A_n \in G$.
Then $v$ extends uniquely to a valuation on $\L$. 
\end{theorem} 

A \emph{(rational) polytope} is the union of finitely many open or closed convex (rational) polytopes. 
%Let $\Rat(d)$ denote the set of all rational polytopes of dimension $\leq d$. Note that $\Rat(d)$ is a lattice of sets. 
We apply Groemer's Theorem as follows. Given a (rational) polytope $\P$, fix a triangulation of $\P$ into (rational) simplices. Let $\T_\P$ denote the set of all unions of faces of these simplices. Hence $\T_\P$ is a lattice, whose top element is $\P$. The generating set of $\T_\P$ consists of the simplices of the triangulation and their faces, which are also simplices. The union of $n$ faces $\F_1, \F_2, \dots, \F_n$ is again a face if only if one of the faces $\F_k$ contains all the others. In this situation the condition of Groemer's Theorem, Theorem~\ref{T:groemer}, is directly satisfied. This means that we can define a function arbitrarily on the faces of the triangulation, and this function will extend uniquely to a valuation on $\T_\P$. Moreover, the values on elements of $\T_\P$ are given by an iterated application of the inclusion-exclusion formula (\ref{inclexcl}). In particular:

\begin{corollary}\label{vollcor}
Two valuations that agree on all (rational) simplices agree on all (rational) polytopes.
\end{corollary}
This result is useful to us, since $V_\P$ is a valuation as a function of $\P$. 

%%%%%%%%%%%%%%%%%%%%%%%%%%%%%%%%%%%%%%%%%%%%%%%%%%%%%%%%%%%%%%%%%%%%%%%%%%%%%%
\newcommand{\w}{{\bf w}}

\section{Proof of the reciprocity theorem}
\setcounter{equation}{0}

\begin{proof}[Proof of Theorem~\ref{T:mainthm}] 
Because of Corollary~\ref{vollcor}, it is enough to prove Theorem~\ref{T:mainthm} for simplices. 
Without loss of generality we may assume that
the simplex $\S$ has the same dimension as the space it lies in.
Hence suppose $\S$ is a rational $d$-simplex in $\R^{d}$, whose vertices have coordinates with denominator $p$. We will prove the reciprocity identity (\ref{E:mainrec}) for the function $V_\S (r+pt)$ for a fixed~$r$. The fact that $V_\S (r+pt)$ is a polynomial in $t$ will be recovered in passing. 

By Theorem~\ref{T:brionext}, 
  \[ V_\S ( r + pt ) \ = \ \sum_{\m \in (r+pt)\S \cap \Z^d} \!\!\! v_{ (r+pt)\S } (\m) \ = \ \lim_{\z \to 1} \beta_{(r+pt)\S} (\z) \ = \ \lim_{\z \to 1} \sum_{ \v \text{ \rm vertex of } \S } \beta_{(r+pt)\S_\v} (\z) \ , \] 
so we need to look at the integer generating functions for $(r+pt)$-dilates of cones more closely. For a vertex $\v$ of $\S$, suppose $\S_\v = \v + \sum_{ k=1 }^{ d } \R_{ \geq 0 } \w_k $ for some integer vectors $\w_1, \dots, \w_d$; then 
  \[ (r+pt) \S_\v \ = \ (r+pt) \v + \sum_{ k=1 }^{ d } \R_{ \geq 0 } \w_k \ = \ tp \v + \left( r \v + \sum_{ k=1 }^{ d } \R_{ \geq 0 } \w_k \right) .\] 
Note that $p \v$ is an integer vector. Now let $R_\v (\z)$ denote the rational function equal to the integer generating function of $r \v + \sum_{ k=1 }^{ d } \R_{ \geq 0 } \w_k$; because this cone is simple, $R_\v (\z)$ is easy to write down: let $\underline \Pi_\v = r \v + \sum_{ k=1 }^{ d } [0,1) \w_k$, then 
  \[ R_\v (\z) = \frac{ \beta_{ \underline \Pi_\v } (\z) }{ \prod_{ k=1 }^{ d } \left( 1 - \z^{ \w_k }  \right)  } \ , \] 
whence
  \[ \beta_{(r+pt) \S_\v} (\z) = \z^{tp \v} R_\v (\z) = \z^{tp \v} \frac{ \beta_{ \underline \Pi_\v } (\z) }{ \prod_{ k=1 }^{ d } \left( 1 - \z^{ \w_k }  \right)  } \ . \] 
Note that we are using the fact that $v$ is a valuation.
For the open cone $\S_\v^\circ$ we obtain, completely analogously,
  \[ \beta_{(r+pt) \S_\v^\circ} (\z) = \z^{tp \v} \frac{ \beta_{ \overline \Pi_\v } (\z) }{ \prod_{ k=1 }^{ d } \left( 1 - \z^{ \w_k }  \right)  } \ , \] 
where $\overline \Pi_\v = r \v + \sum_{ k=1 }^{ d } (0,1] \w_k$.

From the form of these generating functions, we can immediately conclude that $V_\S ( r + pt )$ is a polynomial in $t$ which implies that $V_\S$ is a quasi-polynomial: We know that the sum of the generating functions of all vertex cones is a 
polynomial in the variables of $\z$, hence the singularities of the rational functions cancel. To compute
  \[ V_\S ( r + pt ) \ =\ \lim_{\z \to 1} \sum_{ \v \text{ \rm vertex of } \S } \z^{tp \v} R_\v (\z) \ , \] 
we can write all the rational functions on the right-hand side over one denominator and use L'Hospital's Rule to compute the limit. The result is a polynomial in $t$, as we simply evaluate $\z$ at 1 after using L'Hospital's Rule the correct number of times.

To prove the reciprocity law, we relate the geometry of $\underline \Pi_\v$ with the geometry of $\overline \Pi_\v$. This geometry depends on $r$; let us include this dependency in our notation. Recall that 
\[ \underline \Pi_\v = \underline \Pi_\v (r) = r \v + \sum_{ k=1 }^{ d } [0,1) \w_k \] 
and 
\[ \overline \Pi_\v = \overline \Pi_\v (r) = r \v + \sum_{ k=1 }^{ d } (0,1] \w_k \ . \] 
The two half-open parallelepipeds relate as
\begin{equation}\label{E:parallel}
  \overline \Pi_\v (r) = - \underline \Pi_\v (-r) + \sum_{k=1}^d \w_k \ ,
\end{equation}
as illustrated in Figure~\ref{parallel}.
Using the notation $\z^{-1} = \left( z_1^{-1}, \dots, z_d^{-1} \right)$, we have 
\[ \beta_{ - \underline \Pi_\v (-r) } (\z) = \beta_{ \underline \Pi_\v (-r) } \left( \z^{ -1 }  \right) , \]
which allows us to rephrase (\ref{E:parallel}) as
  \[ \beta_{ \overline \Pi_\v (r) } (\z) = \beta_{ \underline \Pi_\v (-r) } \left( \z^{ -1 }  \right) \prod_{ k=1 }^{ d } \z^{ \w_k } \ . \]

\begin{figure}[htb]
\begin{center}
\includegraphics[totalheight=3.5in]{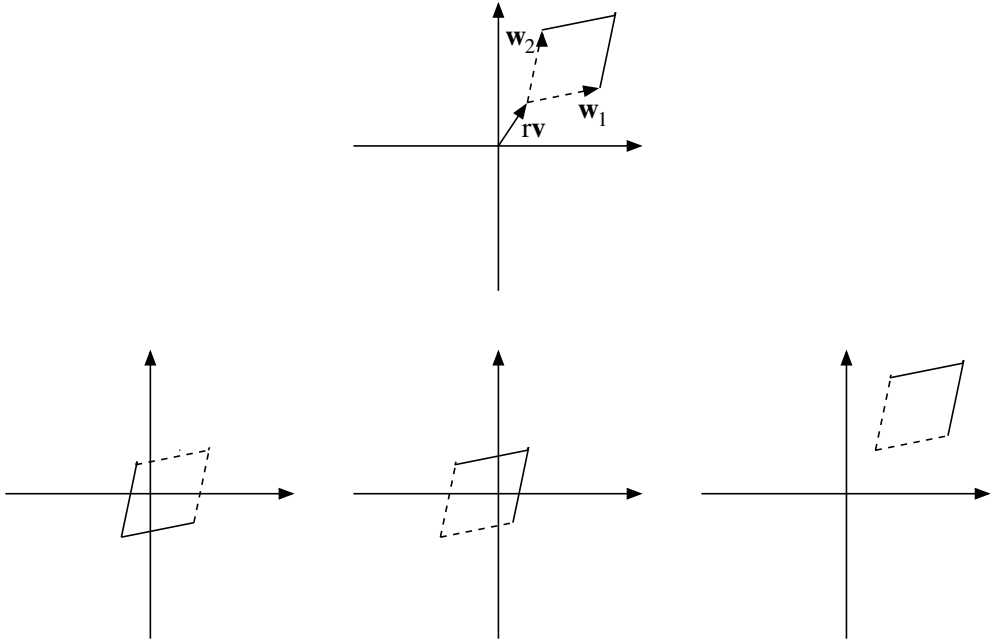}
\end{center}
\caption{Top: $\overline \Pi_\v (r)$ \ Bottom: $\underline \Pi_\v (-r) \ \rightarrow \ -\underline \Pi_\v (-r) \ \rightarrow \ - \underline \Pi_\v (-r) + \sum \w_k$}\label{parallel}
\end{figure} 

Now by the extended Brion-Theorem for simplices (Theorem~\ref{T:brionext}) and the rational generating functions for simple cones, 
\begin{align*} V_{ - \S^\circ }  ( r+pt ) &= \lim_{\z \to 1} \beta_{ - (r+pt) \S^\circ } (\z) \\
  &= \lim_{\z \to 1} \beta_{ (r+pt) \S^\circ } \left( \z^{ -1 } \right) \\
  &= \lim_{\z \to 1} \sum_{ \v \text{ \rm vertex of } \S } \z^{-tp \v} \frac{ \beta_{ \overline \Pi_\v (r) } \left( \z^{ -1 }  \right) }{ \prod_{ k=1 }^{ d } \left( 1 - \z^{ -\w_k }  \right)  } \\
  &= \lim_{\z \to 1} \sum_{ \v \text{ \rm vertex of } \S } \z^{-tp \v} \frac{ \beta_{ \underline \Pi_\v (-r) } \left( \z \right) \prod_{ k=1 }^{ d } \z^{ -\w_k } }{ \prod_{ k=1 }^{ d } \left( 1 - \z^{ -\w_k }  \right)  } \\
  &= \lim_{\z \to 1} \sum_{ \v \text{ \rm vertex of } \S } \z^{-tp \v} \frac{ \beta_{ \underline \Pi_\v (-r) } \left( \z \right) }{ \prod_{ k=1 }^{ d } \left( \z^{ \w_k } - 1 \right)  } \\
  &= (-1)^d V_\S (-r-pt)
\end{align*}
\end{proof}

%%%%%%%%%%%%%%%%%%%%%%%%%%%%%%%%%%%%%%%%%%%%%%%%%%%%%%%%%%%%%%%%%%%%%%%%%%%%%%

\section{Concluding remarks}
\setcounter{equation}{0}

There exists a slightly more general version of our main Theorem~\ref{T:mainthm}. Namely, instead of $\P$ one can take $\P^\circ$ together with some facets of $\P$. Then $\P^\circ$ on the other side of the reciprocity identity gets replaced by $\P^\circ$ together with the remaining facets. This parallels a theorem of Stanley which generalizes Ehrhart-Macdonald Reciprocity~\cite{stanleyreciprocity}.

There is also a theorem analogous to Brion's for polyhedra (not just polytopes). Our proof goes through in this case once one allows \emph{unbounded simplices}, which result when one moves a certain face of a simplex to infinity. (Any polyhedron can be triangulated into bounded and unbounded simplices.)

%%%%%%%%%%%%%%%%%%%%%%%%%%%%%%%%%%%%%%%%%%%%%%%%%%%%%%%%%%%%%%%%%%%%%%%%%%%%%%
\bibliographystyle{amsplain}
%\bibliography{bib}

\def\cprime{$'$}
\providecommand{\bysame}{\leavevmode\hbox to3em{\hrulefill}\thinspace}
\providecommand{\MR}{\relax\ifhmode\unskip\space\fi MR }
\providecommand{\MRhref}[2]{%
  \href{http://www.ams.org/mathscinet-getitem?mr=#1}{#2}
}
\providecommand{\href}[2]{#2}

\setlength{\parskip}{0cm} 
\end{document}